\input amstex
\input amssym
\documentstyle{amsppt}
\def\e{\epsilon}
\def\a{\alpha}
\def\b{\beta}
\def\th{\theta}

\def\s{\sigma}

\def\o{\overline}
\def\f{\flushpar}

\def\om{\omega}
\def\Om{\Omega}
\def\B{\Cal B}

\def\C{\Cal C}

\def\T{\widehat T}
\def\mathcal{\Cal}
\def\({\biggl(}
\def\){\biggr)}

\def\<{\bold\langle}
\def\>{\bold\rangle}
\def\pprime{\prime\prime}

\document\topmatter
\title On the mixing coefficients of piecewise monotonic  maps.\endtitle
\author J. Aaronson,\   H. Nakada \endauthor
 \address[Jon. Aaronson]{\ \ School of Math. Sciences, Tel Aviv University,
69978 Tel Aviv, Israel.}
\endaddress
\email{aaro\@tau.ac.il}\endemail\address[Hitoshi Nakada]{\ \ Dept.
of Math., Keio University,Hiyoshi 3-14-1 Kohoku,
 Yokohama 223, Japan}\endaddress
\email{nakada\@math.keio.ac.jp}\endemail\abstract We investigate
the mixing coefficients of interval maps satisfying Rychlik's
conditions. A mixing Lasota-Yorke map is  reverse $\phi$-mixing.
If its invariant density is uniformly  bounded away from $0$, it
is $\phi$-mixing iff all images of all orders are big in which
case it is $\psi$-mixing. Among $\b$-transformations,
 non-$\phi$-mixing is generic. In this sense, the asymmetry of $\phi$-mixing is natural.

 \endabstract\thanks{\copyright 2004 Preliminary version.}\endthanks\subjclass 37A25,\ 37E05 (37D20, 37C30, 37C40, 60G10)\endsubjclass
\endtopmatter \heading\S0 Introduction\endheading
\subheading{Mixing and measures of dependence between
$\s$-algebras} A mixing property of a stationary stochastic
process $(\dots, X_{-1},X_0,X_1,\dots)$ reflects a decay of the
statistical dependence between the past $\s$-algebra $\s(\{X_k:\
k\le 0\})$ and the  asymptotic future $\s$-algebra $\s(\{X_k:\
k\ge n\})$ as $n\to\infty$ and the various mixing properties are
described by corresponding measures of dependence between
$\s$-algebras (see \cite{Br}).

Let $(\Om,\mathcal F,P)$ be a probability space, and let $\mathcal
A,\mathcal B\subset\mathcal F$ be sub-$\s$-algebras. We consider
the following measures of dependence:
$$\psi(\mathcal A,\mathcal B):=\sup\,\{|\tfrac{|P(A\cap B)-P(A)P(B)}{P(A)P(B)}|:\ B\in\mathcal B_+,\
A\in\mathcal A_+\}$$
$$\phi(\mathcal A,\mathcal B):=\sup\,\{|\tfrac{|P(A\cap B)}{P(A)}-P(B)|:\ B\in\mathcal B_+,\
A\in\mathcal A_+\}$$
$$\align\b(\mathcal A,\mathcal B):=\sup\{ &\tfrac12\sum_{a\in\zeta,\ b\in\xi}|P(a\cap
b)-P(a)P(b)|:\ \\ &\zeta\subset\mathcal A,\ \xi\subset\mathcal B\
\text{\rm finite partitions}\}\endalign$$ where, for
$\C\subset\Cal F,\ \C_+:=\{C\in\C:\ m(C)>0\}$.

 As shown in
\cite{Br},
$$\b(\mathcal A,\mathcal B)\le\phi(\mathcal A,\mathcal B)\le
\tfrac12\psi(\mathcal A,\mathcal B).$$ Note that $\b(\mathcal
A,\mathcal B)=\b(\mathcal B,\mathcal A)$ and $\psi(\mathcal
B,\mathcal A)= \psi(\mathcal A,\mathcal B)$ but  it may be that
$\phi(\mathcal B,\mathcal A)\ne\phi(\mathcal A,\mathcal B)$ (see
\cite{Br} and below).

Accordingly, we set $\phi_+(\mathcal A,\mathcal B):=\phi(\mathcal
A,\mathcal B)$ and $\phi_-(\mathcal A,\mathcal B):=\phi(\mathcal
B,\mathcal A).$

These measures of dependence give rise to the following {\it
mixing coefficients} of a stationary stochastic process $(\dots,
X_{-1},X_0,X_1,\dots)$ defined on the probability space
$(\Om,\mathcal F,P)$:
$$\xi(n):=\xi(\s(\{X_t\}_{ t\le 0}),\s(\{X_t\}_{t\ge n+1}))\ \ \ (\xi=\psi, \phi_+, \phi_-, \b).$$
and the stationary stochastic process $(\dots,
X_{-1},X_0,X_1,\dots)$ is called {\it $\xi$-mixing}  if
$\xi(n)\underset{n\to\infty}\to\rightarrow 0$.

The property $\b$-mixing is also known as {\it absolute
regularity} and {\it weak Bernoulli}.

By stationarity, we have
$$\xi(n)=\sup_{k\ge 1}\xi(\s(\{X_t\}_{0\le t\le k-1}),\s(\{X_t\}_{t\ge n+k+1}))\ \ \ (\xi=\psi, \phi_+, \phi_-, \b).$$

\subheading{ Piecewise monotonic  maps}
A {\it non-singular, piecewise monotonic map } (PM map) of the
interval $X:=[0,1]$  is denoted $(X,T,\a)$ where  $\a$ is a finite
or countable collection of open subintervals of $X$ which is a
partition in the sense that $\bigcup_{a\in\a}a=X\ \mod m$ (where
$m$ is Lebesgue measure) and $T:X\to X$ is a map such that $T|_A$
is absolutely continuous, strictly monotonic
 for each $A\in\a$. \par Let $(X,T,\a)$ be a PM map. For each $n\ge 1$,  $(X,T^n,\a_n)$ is
also a PM map, where
$$\a_n:=\{[a_0,\dots,a_{n-1}]:=\bigcap_{k=0}^{n-1}T^{-k}a_k:\
a_0,\dots,a_{n-1}\in\a\}.$$
\par A PM map $(X,T,\a)$ satisfies
$m\circ T^{-1}\ll m$, whence $f\in L^\infty(m)\ \Rightarrow\
f\circ T\in L^\infty(m)$ . Let $\T:L^1(m)\to L^1(m)$ be the
predual of $f\mapsto f\circ T$, then
$$\T^ng=\sum_{a\in\a_n}v_a'1_{T^na}g\circ v_a$$  where $v_a:T^na\to a$ is given by
$v_a:=(T^n|_a)^{-1}$. Under certain additional assumptions (see
below), $\exists\ h\in L^1(m),\ h\ge 0,\ \int_Xhdm=1,$ so that $\T
h=h$, i.e. $dP=hdm$ is an  absolutely continuous, $T$-invariant
probability (a.c.i.p.). \subheading{ Mixing coefficients of PM
maps}
\par Now let  $(X,T,\a)$ be a PM map with a.c.i.p.
 $P$. We define the  {\it mixing coefficients}  with respect to the probability space $(X,\B(X),P)$:

 $$\xi(n):=\sup_{k\ge 1}\xi(\s(\a_0^{k-1}), T^{-(n+k)}\B(X))\ \ \ (\xi=\psi, \phi_+, \phi_-, \b)$$
and call the PM map {\it $\xi$-mixing}  if $\xi(n)\underset{n\to\infty}\to\rightarrow 0$.

\subheading{Example: Gauss' continued fraction map}
\par This PM map $(X,T,\a)$ with $Tx=\tfrac1x\ \mod 1$ and $h(x)=\tfrac1{\ln 2(1+x)}$ was one of the first considered (see \cite{I-K}).
 Following
Kuzmin's proof (in \cite{Ku}) that  $\|\T^n1-h\|_\infty\le
M\th^{\sqrt n}$ (some $M>0,\ \ \th\in (0,1)$), Khinchine noted  (\cite{Kh}) that
$$\psi^{\star}(n):=\underset{a\in\a_k,\ k\ge
1}\to{\sup}\|\T^n\tfrac{v'_a}{m(a)}-h\|_\infty\le M\th^{\sqrt n}.$$ This estimate was improved by L\' evy to $\psi^{\star}(n)\le M\th^{n}$ (see \cite{L}).

To see  that $(X,T,\a)$ is $\psi$-mixing, one estimates the similar
$$\psi^{\circ}(n):=\underset{a\in\a_k,\ k\ge
1}\to{\sup}\tfrac1{m(a)}\|\T^n({v'_ah\circ v_a})-hP(a)\|_\infty.$$

   The exponential convergence to zero of $\psi^{\circ}(n)$ was shown in \cite{Go}. Theorem 1(b) is a generalization of this
to non-Markov situations. The connection with mixing is seen
through the identity $\T^{n+k}({1_ah})=\T^n({v'_ah\circ v_a})$ for
$k,n\ge 1,\ a\in\a_k$ (see below).

\subheading{RU maps}\par These are PM maps
satisfying the conditions (U) and (R) below. \par The PM map
$(X,T,\a)$ is called {\it uniformly
 expanding} if
 $$\inf_{x\in a\in\a}|T'(x)|=:\b>1;\tag{U}$$
 and is said to \f$\bullet\ \ \ $ satisfy {\it
Rychlik's condition} (\cite{Ry}, see  also \cite{ADSZ}) if
$$
\sum_{A\in\a}\|1_{TA}v_A'\|_{BV}=:\Cal R<\infty.\tag{R}
$$
where $\|f\|_{\text{\rm BV}}:=\|f\|_\infty+\bigvee_Xf$ and
  $v_A'$ is that version of this $L_1$-function
which minimizes variation. \par Suppose that $(X,T,\a)$ satisfies
(R) and (U), then by proposition 1 in \cite{Ry},
 $\exists\ M_0>0,\ \th\in (0,1)$ such that
$$\sum_{a\in\a_n}\|1_{T^na}v_a'f\circ v_a\|_{\text{\rm BV}}\le
M_0(\th^n \|f\|_{BV}+\|f\|_1)\le 2M_0\|f\|_{BV}\ \forall\ n\ge 1.
\tag1$$ It follows (see \cite{Ry}) that the ergodic decomposition
of $(X,\B(X),m,T)$ is finite, that each ergodic component is open
$\mod\ m$ and that the tail $\bigcap_{n=0}^\infty T^{-n}\B$ is
finite and cyclic on each ergodic component. Moreover on each
ergodic component $C,\ \exists\ h=h_C:C\to\Bbb R,\ h\ge 0,\
\int_Xhdm=1,\ h\in BV,\ [h>0]$ open $\mod\ m$ so that $\T h=h$ and
$C=\bigcup_{n=0}^\infty T^{-n}[h>0]\ \mod\ m$.  The probability
$dP_C=h_Cdm$ is an ergodic a.c.i.p. for $T$. \par A  RU map
$(X,T,\a)$ which is conservative and ergodic with respect to $m$
is called {\it basic}. In this case, $\exists\ \ h:X\to\Bbb R,\
h>0,\ \int_Xhdm=1,\ h\in BV$ so that $\T h=h$.
\par An ergodic  RU map
$(X,T,\a)$ (not necessarily basic) is called {\it weakly mixing}
if $(X,\B(X),P,T)$ is weakly mixing (where $P\ll m$ is the
a.c.i.p. for $T$).

Let $(X,T,\a)$ be a weakly mixing RU map, then it is exact with
respect to $m$, and $\exists C>0,\ r\in (0,1)$ so that
$$\|\T^{n}f-\int_Xfdmh\|_{\text{\rm BV}}\le Cr^n\|f\|_{\text{\rm BV}}\
\forall\ f\in\text{\rm BV},\  n\ge 1\tag2$$ (\cite{Ry}, see also
the earlier \cite{H-K} for the case where $\#\a<\infty$).
\subheading{AFU maps}
\par  The PM  map $(X,T,\a)$ is called $C^2$ if,  for all
$A\in\a$, $T:\o{A}\to T\o{A}$ is a $C^2$ diffeomorphism. The
$C^2$  PM $(X,T,\a)$ map is called  an {\it AFU map} (as in
\cite{Z}) if it satisfies (U),
$$T\a:=\{TA:A\in\a\}\ \ \text{\rm is finite};\tag{F}$$ and
$$\sup_X\tfrac{|T''|}{(T')^2}<\infty.\tag{A}$$
A {\it Lasota-Yorke (LY) map} is an AFU map $(X,T,\a)$ with $\a$ finite (as in \cite{L-Y}).

\bigskip

 \par Let $(X,T,\a)$ be an AFU map, then as can be gleaned from  \cite{Z},
 \f$\bullet\ \ \ $  $(X,T,\a)$ is a RU map whose ergodic
 components are finite unions of intervals,
  \f$\bullet\ \ \ $ $1_{[h>0]}\tfrac1h\in \text{\rm BV}$ (and $\tfrac1h\in \text{\rm BV}$ in case
  $(X,T,\a)$ is basic) whenever $h\in BV,\ h\ge  0,\ \T h=h$,
   \f$\bullet\ \ \ $ $\exists
K>0$ so that $$|v_a^{\pprime}(x)|\le Kv_a^{\prime}(x),\
v_a'(x)=e^{\pm K}\tfrac{m(a)}{m(T^na)}\ \forall\ n\ge 1,\
a\in\a_n.\tag3$$
  \heading \S1 Mixing coefficients of RU maps\endheading
As shown in \cite{Ry}, if $(X,T,\a)$ is a weakly mixing RU map,
then  $$\b(n)\le 2CM_0\|h\|_{\text{\rm BV}}r^n\ \ \ (n\ge 1)$$
where $M_0$ is as in (1)  and $C>0,\ r\in (0,1)$
 are as in (2)  (see the
remark after lemma 2 below). \proclaim{Theorem 1}\par {\rm (a)}
Let $(X,T,\a)$ be a weakly mixing RU map. \f If \
$\inf_{[h>0]}\,h>0$, then $\exists\ B>0$ so that $\phi_-(n)\le
Br^n$.\par {\rm (b)} Let $(X,T,\a)$ be a basic, weakly mixing AFU
map. \f If  \ $\inf_{n\ge 1,\ a\in\a_n}{m(T^na)}>0$, then
$\exists\ B>0$ so that $\psi(n)\le Br^n$.\par {\rm (c)} Let
$(X,T,\a)$ be a PM map with a.c.i.p. $P\ll m,\ \#\ \a<\infty$ and
so that $\s(\{T^{-n}\a:\ n\ge 0\})=\B(X)$. \f If\   $\inf_{n\ge
1,\ a\in\a_n}{m(T^na)}=0$, then $\phi_+(n)=1\ \forall\ n\ge 1
$.\endproclaim\subheading{Remarks}\par 1) It  follows that
$\phi_-(n)\to 0$ exponentially for any weakly mixing  AFU map.
This result was announced in \cite{Go} for
$\b$-transformations.\par 2) Part (b) of the theorem is only
established for basic maps as we do not know whether $\inf_{n\ge
1,\ a\in\a_n}{m(T^na)}>0$ implies
$$\inf\,\{m(T^na\cap [h>0]):\ n\ge 1,\ a\in\a_n,\ m(a\cap [h>0])>0
\}>0.$$
 \heading{Proof of theorem 1
}\endheading \proclaim{Lemma 2} Let $(X,T,\a)$ be a weakly mixing
RU map, then
$$|P(A\cap
T^{-(n+k)}B)-P(B)P(A)|\le Mr^nm(B\cap [h>0])$$
$\forall\ n,k\ge 1,
A\in\s(\a_k),\ B\in\B(X) $
 where
$M:=2CM_0\|h\|_{\text{\rm BV}}$ with $M_0$ as in {\rm (1)}  and
$C>0,\ r\in (0,1)$  as in {\rm (2)}.\endproclaim
 \demo{Proof}  We show first that
 $$\phi_-^\circ(n):=\sup\{\|\T^{n+k}(1_Ah)-P(A)h\|_\infty:\ k\ge 1,\ A\in\s(\a_k)\}
 \le Mr^n\ \forall\ n\ge 1.\tag4$$ The   sequence $\phi_-^\circ(n)$ is the analogue of
$\psi^\circ(n)$ for $\phi_-$-mixing.
\par To see (4), fix $k\ge 1$ and suppose that $A\in\s(\a_k),\ A=\bigcup_{a\in\goth a}a$ where
$\goth a\subseteq\a_k$, then,
$$\T^k(h1_A)=\sum_{a\in\goth a}v_a'1_{T^ka}h\circ v_a,$$
$$P(A)=\sum_{a\in\goth a}P(a)=\sum_{a\in\goth
a}\int_{T^ka}v_a'h\circ v_adm,$$ and for $n\ge 1$,

$$\T^{n+k}(1_Ah)=\sum_{a\in\goth a}\T^n(v_a'1_{T^ka}h\circ v_a).$$
 Thus
$$\align \|\T^{n+k}(1_Ah)-P(A)h\|_\infty &\le
\sum_{a\in\goth
a}\|\T^n(v_a'1_{T^ka}h\circ v_a)-h\int_{T^ka}v_a'h\circ v_adm\|_\infty\\
& \le Cr^n\sum_{a\in\goth a}\|v_a'1_{T^ka}h\circ v_a\|_{\text{\rm
BV}}\ \ \text{\rm by (2)}\\ &\le 2CM_0\|h\|_{\text{\rm BV}}r^n\ \
\text{\rm by (1)}\endalign$$ establishing (4). Using (4), for
$k\ge 1,\ A\in\a_k,\ B\in\B$,
$$\align |P(A\cap T^{-(n+k)}B)-P(A)P(B)| &=|\int_{B\cap [h>0]}(\T^{n+k}(1_Ah)-P(A)h)dm|\\ &\le Mr^nm(B\cap [h>0]).\endalign$$
\hfill\qed\enddemo
\subheading{Remark}\ \ It was shown in \cite{Ry}, using a version of lemma 2, that for
 a weakly
mixing RU map,  $\b(n)\le Mr^n$.  \demo{Proof of  theorem 1 {\rm
(a)}} By lemma 2, for $n,k\ge 1, A\in\s(\a_k),\ B\in\B(X)$,
$$|P(A\cap
T^{-(n+k)}B)-P(B)P(A)|\le Mr^nm(B\cap [h>0])\le M\|\tfrac1h\|_{L^\infty([h>0])}
r^nP(B).$$\hfill\qed\enddemo
 \demo{Proof of  theorem 1 {\rm (b)}}

 We
show first that if $\inf_{n\ge 1,\ a\in\a_n}{m(T^na)}=:\eta>0$,
then
 $$\|\T^{n+N}(1_ah)-P(a)h\|_\infty\le M_1r^nm(a)\ \forall\ N,n\ge 1,\
a\in\a_N\tag5$$ where $M_1:=\tfrac{4Ce^K\|h\|_{\text{\rm
BV}}}{\eta}$. A standard calculation shows that
$$\|\T^k(1_ah)\|_{\text{\rm
BV}}=\|v_a'1_{T^ka}h\circ
v_a\|_{\text{\rm BV}}\le 4e^K\|h\|_{\text{\rm
BV}}\tfrac{m(a)}{m(T^ka)}\ \forall\ k\ge 1,\ a\in\a_k.\tag6$$ To
see (5), fix $N,n\ge 1,\ a\in\a_N$, and note that
$\T^{N}(1_ah)=1_{T^Na}h\circ v_a v_a'$. By (2),(3) and (6),
$$\align \|\T^{n+N}(1_ah)-P(a)h\|_\infty &=
\|\T^{n}(1_{T^Na}h\circ v_a v_a')-P(a)h\|_\infty\\ &\le
Cr^n\|1_{T^Na}h\circ v_a v_a'\|_{\text{\rm BV}}\\ &\le
4Ce^K\|h\|_{\text{\rm BV}}\tfrac{m(a)}{m(T^Na)}r^n \le M_1r^n
m(a).\endalign
$$ Now let $N\ge 1$ and suppose that $A\in\s(\a_N),\ A=\bigcup_{a\in\goth a}a$ where
$\goth a\subseteq\a_N$. It follows from (5) that $\forall\
B\in\B,\ n\ge 1$, $$\align |P(A\cap T^{-(n+N)}B)-P(A)P(B)| &\le
\int_{B}|\T^{n+N}(1_Ah)-P(A)h|dm\\ &\le m(B)\|\T^{n+N}(1_Ah)-P(A)h
\|_\infty\\ & \le M_1r^nm(B)\sum_{a\in\goth a}m(a)\\
&=M_1r^nm(A)m(B)\\ &\le  M_1(\|\tfrac1h\|_{\infty})^2
r^nP(A)P(B).\qed\endalign$$ \enddemo

 \demo{Proof of theorem 1 {\rm
(c)}}
 \par  Since $\s(\{T^{-n}\a:\ n\ge 0\})=\B(X)$, $\max\{m(a):\ a\in\a_n\}\to 0$ as $n\to\infty$. Fix $0<\e<1$ and choose $\ell\ge 1$ so that
 $\max\{m(a):\ a\in\a_\ell\}<\tfrac{\e}{2\|h\|_\infty}$.

 We show that $\phi_+(N)\ge 1-\e\ \forall\
 N\ge 1$. Indeed, fix $N\ge 1$.

 Since $\inf_{n\ge
1,\ a\in\a_n}{m(T^na)}=0$, and $\#\a_{N+\ell}<\infty,$ $\exists\ k\ge 1,\ \om\in\a_k$ so that
$m(T^k\om)<\min\{m(a):\ a\in\a_{N+\ell}\}$.

Since $T^k\om$ is an interval, $\exists\
b=[b_1,\dots,b_{N+\ell}],\ c=[c_1,\dots,c_{N+\ell}]\in\a_{N+\ell}$
so that $T^k\om\subset J:=b\uplus c$. Next, $T^NJ\subset
[b_{N+1},\dots,b_{N+\ell}]\cup
[c_{N+1},\dots,c_{N+\ell}]\in\s(\a_{\ell})$, whence $m(T^NJ)\le
2\max\{m(a):\ a\in\a_\ell\} <\tfrac{\e}{\|h\|_\infty}$ and
$\exists\ B\in\s(\a_\ell),\ B\cap T^NJ=\emptyset,\
m(B)>1-\tfrac{\e}{\|h\|_\infty}$. It follows that $P(B)>1-\e$, and
that
$$\phi_+(N)\ge P(B)-P(T^{-(N+k)}B|\om)=P(B)>1-\e.\qed$$\enddemo

\heading{\S2 Examples}\endheading \subheading{$\beta$-expansions }
\par For $\b>1$, consider $(X,T_\b,\a)$ where $X=[0,1]$, $T_\b:X\to X$ is given
by $T_\b x:=\{\b x\}$. and
$\a:=\{[\tfrac{j}\b,\tfrac{j+1}\b)\}_{j=0}^{[\b]-1}\cup\{[\tfrac{[\b]}\b,1)\}$.
As shown in \cite{Pa},  $(X,T_\b,\a)$ is a basic LY map. Theorem 1(a) applies and
$\phi_-(n)\to 0$ exponentially. To apply theorem 1(b), we prove:
\proclaim{Proposition 3}
$$\inf_{n\ge 1,\ a\in\a_n}\,m(T_\b^na)=\inf_{n\ge
1}\,T_\b^n1.\tag7$$\endproclaim\demo{Proof}
 \par Define (as in \cite{Pa})
$\pi:X\to \{0,1,\dots,[\b]\}^\Bbb N$ by $\pi(x)_n:=[\b T_\b
^{n-1}x]$ and set $X_\b:=\o{\pi_\b(I)}$. Let $\prec$ denote
lexicographic order on $\{0,1,\dots,[\b]\}^\Bbb N$, then (see
\cite{Pa})  $x< y$ implies $\pi(x)\prec\pi(y)$  and that
$\pi(x)\prec\pi(y)$ implies $x\le y$.

 Let
$$\om=\om_\b:=\cases & (\o{a_1,a_2,\dots,a_{q-1},a_q-1})\ \ \ \
\pi_\b(1)=(a_1,a_2,\dots,a_{q-1},a_q,\o 0),\\ & \pi_\b(1)\ \ \ \
\text{else.}\endcases$$
\par By \cite{Pa},
$$X_\b=\{y\in\{0,1,\dots,[\b]\}^\Bbb N:\ y_k^\infty\prec \om\
\forall\ k\ge 1\}.\tag8$$ where $y_k^\infty:=(y_k,y_{k+1},\dots)$.
For $a=[a_1,\dots,a_N]\in\a_N$, define $$K_N(a):=\cases &0\ \ \ \
\ \ \ \ \ \nexists\ 1\le n\le N,\ a_{N-n+1}^N=\om_1^{n},\\ &
\max\{1\le n\le N,\ a_{N-n+1}^N=\om_1^{n}\}\ \ \ \ \ \ \ \ \ \
\text{\rm else}\endcases$$ (where
$a_j^k:=(a_j,a_{j+1},\dots,a_k)$), then by (8),
$\pi([a_1,\dots,a_N])=\{x\in X_\b:\
 x\prec\om_{K_N(a)+1}^\infty\}=\pi[0,T_\b^{K_N(a)}1)$, whence
$T_\b^N[a_1,\dots,a_N]=[0,T_\b^{K_N(a)}1).$ The proposition
follows from
 this.\hfill\qed\enddemo
 Thus by theorem 1(b), $T_\b$ is $\psi$-mixing iff $\inf_{n\ge
1}\,T_\b^n1>0$, or equivalently (see \cite{Bl}) $X_\b$ is {\it
specified} in the sense that $\exists\ L\ge 1$ so that
$$m(a\cap T_\b^{-(i+L)}b)>0\ \ \forall\
a\in \a_i,\ b\in\a_j.$$ As shown in \cite{S}, the set
$\{\b>1:X_\b\ \text{\rm specified}\}$ is a meagre set of Lebesgue
measure zero and Hausdorff dimension $1$ in $\Bbb R$ and so
exponential $\psi$-mixing occurs for many $\b>1$ for which $X_\b$
is not sofic .

\subheading{"Japanese" continued fractions}
\par Fix $\a\in (0,1]$ and define
$T=T_\a:[\a-1,\a)\circlearrowleft$ by $$T_\a
x:=|\tfrac1{x}|-[|\tfrac1{x}|+1-\a].$$ These maps generalize the
well known  Gauss map $T_1$. For $\a\in (0,1)$, $T_\a$ is a
topologically mixing, basic AFU map; whence by theorem 1,
\f$\bullet\ \ \ $ exponentially reverse $\phi$-mixing, and
\f$\bullet\ \ \ $   exponentially $\psi$-mixing when
$\inf\,\{m(T^na):\ n\ge 1,\ a\in\a_n\}>0$.

Theorem 1(c) does not apply since $\#\a=\infty$. However, as shown
in \cite{N-N}, for Lebesgue a.e. $\tfrac12<\a< 1$, $T_\a$ is not
$\phi_+$-mixing.

\bigskip
 \centerline{\bf References}
 \medskip
\noindent [ADSZ] J. Aaronson, M. Denker, O.  Sarig, R.
Zweim\"uller, Aperiodicity of cocycles and conditional local limit
theorems, {\it Stochastics and Dynamics}, {\bf 4}, (2004), No. 1
31-62

\medskip
\noindent [Bl]\ \  F. Blanchard, $\beta$-expansions and symbolic
dynamics, {\it Theoret. Comput. Sci.} {\bf 65},  (1989), no. 2,
131--141.

\medskip
\noindent [Br]\ \ R. C. Bradley,  Basic properties of strong
mixing conditions. Dependence in probability and statistics
(Oberwolfach, 1985), 165--192, {\it Progr. Probab. Statist.},{\bf
11}, Birkh\"auser Boston, Boston, MA, 1986.

\medskip
\noindent [Go]\ \  M. I.  Gordin, Stochastic processes generated
by number-theoretic endomorphisms (English) {\it Sov. Math.,
Dokl.} {\bf 9}, (1968) 1234-1237 ; (translation from Russian: {\it
Dokl. Akad. Nauk SSSR} {\bf 182}, (1968) 1004-1006).

\medskip
\noindent [H-K]\ \   F. Hofbauer, G. Keller,  Ergodic properties
of invariant measures for piecewise monotonic transformations.
{\it Math. Z. } {\bf 180 }(1982), no. 1, 119--140.

\medskip
\noindent [I-K]\ \
M. Iosifescu, C. Kraaikamp, {\it Metrical theory of continued fractions}, (Mathematics and its Applications, {\bf 547}), Kluwer Academic Publishers, Dordrecht, 2002.

\medskip
\noindent [Kh]\ \ A. Ya.   Khintchine, {\it Continued fractions},
 P. Noordhoff, Ltd., Groningen 1963;
(translation from Russian: : {\it Cepnye drobi}; Gosudarstv.
Izdat. Tehn.-Teor. Lit., Moscow-Leningrad, 1935).

\medskip
\noindent [Ku]\ \ R. Kuzmin,
Sur une probl\'eme de Gauss.
{\it Atti Congresso Bologna}{\bf 6}, 83-89.

\medskip
\noindent [L-Y]\ \  A.  Lasota, J. Yorke,  On the existence of
invariant measures for piecewise monotonic transformations. {\it
Trans. Amer. Math. Soc.} {\bf 186}, 481--488 (1973).

\medskip
\noindent [L]\ \ P. L\'evy, {\it Th\'eorie de l'addition des
variables al\' eatoires. (French) } (Monogr. d. probabilit. calcul
d. probabilit. et ses appl. 1) Paris: Gauthier-Villars. XVII,
(1937).

\medskip
\noindent [N-N]\ \ H. Nakada, R.  Natsui, Some strong mixing
properties of a sequence of random variables arising from
$\alpha$-continued fractions.{\it  Stoch. Dyn.}{\bf 3}, (2003),
no. 4, 463--476.

\medskip
\noindent [Pa]\ \  W. Parry, On the $\beta $-expansions of real
numbers,{\it Acta Math. Acad. Sci. Hungar.}{\bf 11} (1960)
401--416.

\medskip
\noindent [Ry] M. Rychlik: {Bounded variation and invariant
measures.} {\it Studia Math.} {\bf 76} (1983) 69-80.

\medskip
\noindent [S] J. Schmeling,  Symbolic dynamics for $\beta$-shifts
and self-normal numbers. {\it Ergodic Theory Dynam. Systems} {\bf
17} (1997), no. 3, 675--694.

\medskip
\noindent [Z] R. Zweim\"uller:  Ergodic structure and invariant
densities of non-Markovian interval maps with indifferent fixed
points. {\it Nonlinearity }{\bf 11} (1998), no. 5, 1263--1276.

\enddocument